\documentclass[11pt]{article}

\usepackage[truedimen,margin=25truemm]{geometry}

\usepackage{additional_def}

\usepackage{graphicx}

\usepackage{amsmath}
\usepackage{amsthm}
\newtheorem{Def}{Definition}[section]
\newtheorem{Thm}[Def]{Theorem}

\title{Stochastic variance reduced multiplicative update for \\nonnegative matrix factorization}
\date{\today}

\author{Hiroyuki Kasai\thanks{Graduate School of Informatics and Engineering, The University of Electro-Communications, Tokyo, Japan ({\tt kasai@is.uec.ac.jp}).} }

\begin{document}

\maketitle
\begin{abstract}
Nonnegative matrix factorization (NMF), a dimensionality reduction and factor analysis method, is a special case in which factor matrices have low-rank nonnegative constraints. Considering the stochastic learning in NMF, we specifically address the multiplicative update (MU) rule, which is the most popular, but which has slow convergence property.  This present paper introduces on the stochastic MU rule a variance-reduced technique of stochastic gradient. Numerical comparisons suggest that our proposed algorithms robustly outperform state-of-the-art algorithms across different synthetic and real-world datasets. 
\end{abstract}
\section{Introduction}
\label{Sec:Introduction}

Superior performance of nonnegative matrix factorization (NMF) has been achieved in many technical fields. NMF approximates a nonnegative matrix \mat{V} as a product of two nonnegative matrices \mat{W} and \mat{H}. Given $\mat{V}\in \mathbb{R}^{F \times N}_{+}$, NMF requires factorization of the form $\mat{V} \approx \mat{WH}$, where $\mat{W}\in \mathbb{R}^{F \times K}_{+}$ and $\mat{H}\in \mathbb{R}^{K \times N}_{+}$ are nonnegative {\it factor matrices}. $K$ is usually chosen such that $K \ll \min\{F, N\}$, that is, \mat{V} is approximated in the two {\it low-rank} matrices. This problem is formulated as a constrained minimization problem in terms of the Euclidean distance  as 
\begin{eqnarray}
	\label{Eq:problem}
	\min_{\scriptsize \mat{W}, \mat{H}} &&\frac{1}{2} \| \mat{V} - \mat{WH}\|_F^2=\frac{1}{N} \sum_{n=1}^N \frac{1}{2}\| \vec{v}_n - \mat{W} \vec{h}_n\|^2_2,\nonumber \\
	 {\rm s.t.}&& [\mat{W}]_{f,k} \geq 0, [\mat{H}]_{k,n} \geq 0, \ \ \ \forall f,n,k,
\end{eqnarray}
where $\mat{V}=[\vec{v}_1,\ldots, \vec{v}_N]$ and $\mat{H}=[\vec{h}_1,\ldots, \vec{h}_N]$. $[\mat{A}]_{i,j}$ is the $(i,j)$-th entry of \mat{A}. The non-negativity of \mat{V} enables us to interpret the meanings of the obtained matrices \mat{W} and \mat{H} well. This interpretation produces a broad range of applications in machine learning and signal processing such as text mining, image processing, and data clustering, to name a few. However, because problem (\ref{Eq:problem}) is a {\it non-convex} optimization problem, finding its global minimum is NP-hard. For this problem, Lee and Seung proposed a simple but effective calculation algorithm \cite{Lee2001} as \\
\begin{eqnarray}
	\label{Eq:MU_Batch}
	\mat{H} \leftarrow \mat{H} \odot \frac{{\mat{W}}^{T}\mat{V}}{\mat{W}^T\mat{W}\mat{H}}, \ \ \ 
	\mat{W} \leftarrow \mat{W} \odot \frac{\mat{V}{\mat{H}}^{T}}{\mat{W}\mat{H}{\mat{H}}^{T}}, 
\end{eqnarray}
where $\odot$ (resp.$\frac{\  \cdot\ }{\cdot}$) denotes the component-wise product (resp. division) of matrices, which finds a {\it local optimal} solution of (\ref{Eq:problem}). This rule is designated as the {\it multiplicative update} (MU) rule because a new estimate is represented as the product of a current estimate and some factor. The {\it global convergence to a stationary point} is guaranteed under slightly modified update rules or constraints \cite{Lin_IEEENN_2007,Hibi_ICONIP_2011}. Nevertheless, many efficient algorithms have been developed because the MU rule is accompanied by slow convergence  \cite{Lin_NC_2007,Cichocki_IEICETranA_2009,Gillis_NN_2012,Kim_JGO_2014}. Furthermore, considering big data, an {\it online learning} algorithm is preferred in terms of the computational burden and the memory consumption. Designating the former algorithms as {\it batch-NMF}, this {\it online-NMF} has been  investigated actively in several studies \cite{Bucak_PR_2009,Fvotte_NC_2009,Guan_IEEENNLS_2012,Zhao_AISTATS_2017}. Its robust variant has also been assessed \cite{Zhao_ICASSP_2016}. However, they still exhibit slow convergence. {\it Stochastic gradient descent} (SGD) \cite{Robbins_MathStat_1951} has become the method of choice for solving big data optimization problems. Although it is beneficial because of the low and constant cost per iteration independent of $N$, the {\it convergence rate} of SGD is also slower than that of full GD even for the strongly convex case. For this issue, various {\it variance reduction} (VR) approaches that have been proposed recently have achieved superior convergence rates in convex and non-convex functions. 

This paper presents a proposal of a novel stochastic multiplicative update with the VR technique: SVRMU. The paper also explains extension of SVRMU to the accelerated variant (SVRMU-ACC), and the robust variant (R-SVRMU) for outliers. The paper is organized as follows. Section 2 presents details of the variance reduction algorithm in stochastic gradient. Section 3 presents the proposed stochastic variance reduced multiplicative update (SVRMU). Section 4 provides a convergence analysis. Two extensions are detailed in Section 5. In Section 6, exhaustive comparisons suggest that our proposed SVRMU algorithms robustly outperform state-of-the-art algorithms across different synthetic and real-world datasets. It is noteworthy that the discussion presented here is applicable to other distance functions than the Euclidean distance. The Matlab codes are available at \url{https://github.com/hiroyuki-kasai}.

\section{Variance reduction algorithm in stochastic gradient}

An algorithm designated to solve the problem (\ref{Eq:problem}) without nonnegative constraints is begin eagerly sought in the machine learning field. When designating \mat{W} and $\mat{H}$ as $w$, and designating the rightmost inner term of the cost function (\ref{Eq:problem}) as $f_i(w)$, respectively, {\it full gradient decent} (GD) with a stepsize $\eta$ is the most straightforward approach as $w_{t+1}=w_t - \eta \nabla f(w_t)$, where $\nabla f(w_t)$ corresponds to the gradient $g_t$. However, this is expensive especially when $N$ is extremely large. A popular and effective alternative is a stochastic gradient by which $g_t$ is set to $\nabla f_{n_t}(w_t)$ for $n_t$-th ($n_t \in [N]$) sample that is selected uniformly at random, which is called {\it stochastic gradient descent} (SGD). It updates $w_t$ as $w_{t+1}=w_t - \eta \nabla f_{n_t}(w_t)$, and assumes an {\it unbiased estimator} of the full gradient as $\mathbb{E}_{n_t}[\nabla f_{n_t}(w_t)] = \nabla f (w_t)$. Apparently, the calculation cost per iteration is independent of $N$. {\it Mini-batch} SGD uses $g_t=1/|\mathcal{S}_t| \sum_{n_t \in \mathcal{S}_t} \nabla f_{n_t}(w_t)$, where $\mathcal{S}_t$ is the set of samples of size $|\mathcal{S}_t|$. However, because SGD requires a {\it diminishing} stepsize algorithm to guarantee the convergence, SGD suffers from a slow convergence rate. 

To accelerate this rate, the variance reduction (VR) techniques \cite{Johnson_NIPS_2013_s,Roux_NIPS_2012_s,Shalev_JMLR_2013_s,Defazio_NIPS_2014_s,Zhang_SIAMJO_2014_s,Nguyen_ICML_2017} explicitly or implicitly exploit a full gradient estimation to reduce the variance of noisy stochastic gradient, leading to  superior convergence properties. We can regard this approach as a hybrid algorithm of GD and SGD. A representative research among them is Stochastic Variance Reduced Gradient (SVRG) \cite{Johnson_NIPS_2013}. SVRG first keeps $\tilde{w}=w_t^{s-1}$ indexed by $t=0, \cdots, m_{s}-1$ at the end of $(s-1)$-th epoch with $m_{s-1}$ inner iterations. It also sets the initial value of the inner loop in $s$-th epoch as $w_0^s=\tilde{w}$. Computing a full gradient $\nabla f(\tilde{w})$, it randomly selects $n_t^s \in [N]$  for each $\{t,s\} \geq 0$, and computes a {\it modified stochastic gradient} $g_t^s$ as
\begin{eqnarray}
\label{Eq:SVRG}
g_t^s & = & \nabla f_{n_t^s} (w_{t}^{s}) -  \nabla f_{n_t^s} (\tilde{w}^{s})  + \nabla f(\tilde{w}^{s}). 
\end{eqnarray}
For smooth and strongly convex functions, this method enjoys a linear convergence rate as SDCA, SAG and SAGA. 

\section{Stochastic variance reduced multiplicative update (SVRMU)}
\label{Sec:SVRMU}

\begin{algorithm}
\caption{Stochastic variance reduction multiplicative update (SVRMU)}
\label{Alg:SVRMU}
\begin{algorithmic}[1]
\REQUIRE{\mat{V}, maximum inner iteration $m_s>0$.}
\STATE{Initialize $\tilde{{\mat{W}}}^0$ and $\tilde{{\mat{H}}}^0$.}
\FOR{$s=0,1,\ldots$} 
\STATE{Calculate the components of the full gradient $\tilde{\mat{W}}^s\tilde{\mat{H}}^s(\tilde{\mat{H}}^s)^{T}/N$ and $\mat{V} (\tilde{\mat{H}}^s)^{T}/N$}.
\STATE{Store $\mat{W}_0^s = \tilde{\mat{W}}^{s}$.}
	\FOR{$t=0,1,\ldots, m_s-1$} 
	\STATE{Choose $k=n_t^s \in [N]$ uniformly at random.}
	\STATE{Update $\vec{h}_{k}=\vec{h}_{k} \odot ({(\mat{W}^s_{t})}^T \vec{v}_k)/({(\mat{W}^s_{t})}^T \mat{W}^s_{t} \vec{h}_{k})$.}	
	\STATE{Calculate $\mat{Q}_t^s=\mat{W}^s_t\vec{h}_{k}\vec{h}_{k}^{T}  +\vec{v}_k \tilde{\vec{h}}_k^{T} +\tilde{\mat{W}}^s\tilde{\mat{H}}^s (\tilde{\mat{H}}^s)^{T}/N$.}
	\STATE{Calculate $\mat{P}_t^s=\vec{v}_k\vec{h}_{k}^{T} +\tilde{\mat{W}}^s\tilde{\vec{h}}^s_{k}(\tilde{\vec{h}}_{k}^s)^{T}+\mat{V} (\tilde{\mat{H}}^s)^{T}/N$.}	
	\STATE{Calculate the stepsize ratio $\alpha_t^s$.}
	\STATE{Update $\mat{W}^s_{t+1}=\mat{W}^s_{t}-\alpha_t^s \mat{W}^s_{t}/\mat{Q}_t^s \odot (\mat{Q}_t^s - \mat{P}_t^s)$.}	
	\ENDFOR
	\STATE{Set $\tilde{\mat{W}}^{s}=\mat{W}^s_{m_s}$ and $\tilde{\mat{H}}^{s}=\mat{H}$.}
\ENDFOR
\end{algorithmic}
\end{algorithm}

This section first describes the {\it stochastic multiplicative update}, designated as SMU. Then it details the proposed stochastic variance-reduced MU algorithm, i.e., SVRMU. 
The problem setting is the following: we assume that $n_t$-th ($n_t\in [N]$) column of \mat{V}, i.e. $\vec{h}_{n_t}$, is selected at $t$-th iteration uniformly at random. $\vec{h}_{n_t}$ and $\mat{W}$ are updated alternatively by extending (\ref{Eq:MU_Batch}) as 
\begin{eqnarray}
	\label{Eq:MU_Online}	
	\vec{h}_{n_t} \leftarrow \vec{h}_{n_t} \odot \frac{\mat{W}^T \vec{v}_{n_t}}{\mat{W}^T \mat{W} \vec{h}_{n_t}}, \ \ \ \ 
	\mat{W}  \leftarrow  \mat{W} \odot \frac{\vec{v}_{n_t}\vec{h}_{n_t}^{T}}{\mat{W}\vec{h}_{n_t}\vec{h}_{n_t}^{T}}.
\end{eqnarray}

Especially, the MU rule of \mat{W} is regarded as a special case of SGD with an adaptive stepsize of {\it matrix form} of $\mat{S}_t=\alpha\mat{W}/(\mat{W}\vec{h}_{n_t}\vec{h}_{n_t}^{T}) \in \mathbb{R}_+^{F\times K}$ as 
\begin{eqnarray}
	\label{Eq:MU_Online_GD}
	\mat{W} \leftarrow  \mat{W} - \mat{S}_t \odot (\mat{W}\vec{h}_{n_t}\vec{h}_{n_t}^{T} - \vec{v}_{n_t}\vec{h}_{n_t}^{T}), \nonumber
\end{eqnarray}
where $0< \alpha \leq 1$ is the {\it stepsize ratio} that ensures that $\mat{W}$ and $\mat{H}$ are nonnegative when those initial values are nonnegative. The case of $\alpha=1$ produces (\ref{Eq:MU_Online}) exactly.

According to this interpretation, we consider the VR algorithm for SMU. Similarly as SVRG,  SVRMU has a double loop structure. 
By keeping $\tilde{\mat{W}}^s=\mat{W}_t^{s}$ and $\tilde{\mat{H}}^s=\mat{H}$ indexed by $t=0, \cdots, m_{s}-1$ at the end of $(s$-$1)$-th outer loop with $m_{s-1}$ inner iterations, and also by setting the initial value of the inner loop in $s$-th outer loop as $\mat{W}_0^s=\tilde{\mat{W}}^s$,  we compute the components of the full gradient $\tilde{\mat{W}}^s\tilde{\mat{H}}^s(\tilde{\mat{H}}^s)^{T}/N$ and $\mat{V} (\tilde{\mat{H}}^s)^{T}/N$. For each $s \geq 0$ and $t \geq 0$, we first randomly select $n_t^s \in [N]$ and update $\vec{h}_{n_t^s}$ as in (\ref{Eq:MU_Online}). Hereinafter, $k$ is used instead of $n_t^s$ for notation simplicity. Then, we update $\mat{W}^s_{t}$ with an appropriate stepsize $\mat{S}_t^s$ as shown below.
\begin{eqnarray}
	\label{Eq:SVRMU_1_W}
	\mat{W}^s_{t+1} 
	 &=& \mat{W}^s_{t} - \mat{S}_t^s \odot \biggl[
	(\mat{W}^s_t\vec{h}_{k}\vec{h}_{k}^{T} - \vec{v}_{k}\vec{h}_{k}^{T}) \nonumber \\
	&& - \ (\tilde{\mat{W}} \tilde{\vec{h}}^s_{k}(\tilde{\vec{h}}_{k}^s)^{T}    - \vec{v}_{k} (\tilde{\vec{h}}_{t}^s)^{T})  +\frac{\tilde{\mat{W}}  \tilde{\mat{H}}^s  (\tilde{\mat{H}}^s)^{T}   -  \mat{V} (\tilde{\mat{H}}^s)^{T}}{N} 		
	\biggr] \nonumber \\
	        &        =         &\mat{W}^s_{t} - \mat{S}_t^s \odot  \biggl[
	\left(\mat{W}^s_t\vec{h}_{k}\vec{h}_{k}^{T}  +  \vec{v}_{k} (\tilde{\vec{h}}_{k}^s)^{T}   +  \frac{\tilde{\mat{W}}\tilde{\mat{H}}^s(\tilde{\mat{H}}^s)^{T}}{N}\right) \nonumber \\
	&&- \left(\vec{v}_{k}\vec{h}_{k}^{T} + \tilde{\mat{W}} \tilde{\vec{h}}^s_{k}(\tilde{\vec{h}}_{k}^s)^{T}  +\frac{\mat{V} (\tilde{\mat{H}}^s)^{T}}{N}\right)
	\biggr],
\end{eqnarray}
where $\tilde{\mat{H}}^s=[\tilde{\vec{h}}^s_1,\ldots,\tilde{\vec{h}}^s_N]$. Here, we denote $\mat{Q}_t^s \in \mathbb{R}_+^{F\times K}$ as
\begin{eqnarray*}
\mat{Q}_t^s &=& \mat{W}^s_t\vec{h}_{k}\vec{h}_{k}^{T} +  \vec{v}_{k} (\tilde{\vec{h}}_{t}^s)^{T} + \frac{\tilde{\mat{W}} \tilde{\mat{H}}^s  (\tilde{\mat{H}}^s)^{T}}{N}.
\end{eqnarray*}
We also denote $\mat{P}_t^s \in \mathbb{R}_+^{F\times K}$ as 
 \begin{eqnarray*}
\mat{P}_t^s  &=& \vec{v}_{k}\vec{h}_{k}^{T} + \tilde{\mat{W}} \tilde{\vec{h}}^s_{t}(\tilde{\vec{h}}_{t}^s)^{T}+
\frac{\mat{V} (\tilde{\mat{H}}^s)^{T}}{N}.
\end{eqnarray*}
When $\mat{S}_t^s = \alpha_t^s \mat{P}_t^s/\mat{Q}_t^s$ with the stepsize ratio $\alpha_t^s$, the update rule in (\ref{Eq:SVRMU_1_W}) is reformulated as presented below.
\begin{eqnarray}
	\label{Eq:SVRMU_2_W}
	\mat{W}^s_{t+1} &= &  \mat{W}^s_{t}-\frac{\alpha \mat{W}^s_{t}}{\mat{Q}_t^s} \odot (\mat{Q}_t^s - \mat{P}_t^s).
\end{eqnarray}
The overall algorithm is summarized in Algorithm \ref{Alg:SVRMU}. Additionally, the straightforward extension to the mini-batch variant of Algorithm \ref{Alg:SVRMU} is  defined as 
\begin{eqnarray*}
	\label{Eq:MiniBatch_SVRMU_1_W}
	\mat{W}^s_{t+1}    &      =      &   
 \mat{W}^s_{t} - \mat{S}_t^s   \odot   \biggl[ 
	\left(\frac{\mat{W}^s_t\vec{h}_{k}\vec{h}_{k}^{T}    +  \vec{v}_{k} (\tilde{\vec{h}}_{k}^s)^{T}}{b}  +  \frac{\tilde{\mat{W}}\tilde{\mat{H}}^s(\tilde{\mat{H}}^s)^{T}}{N}   \right) \nonumber \\
	&&- \left( \frac{\vec{v}_{k}{\vec{h}_{k}}^{T} + \tilde{\mat{W}} \tilde{\vec{h}}^s_{k}(\tilde{\vec{h}}_{k}^s)^{T}}{b}+\frac{\mat{V} (\tilde{\mat{H}}^s)^{T}}{N}\right)
	\biggr],	       
\end{eqnarray*}
where $b\ (\leq N)$ is the mini-batch size. $\mat{Q}_t^s$ and $\mat{P}_t^s$ in (\ref{Eq:SVRMU_2_W}) are modified accordingly. 

\section{Convergence analysis}
The convergence analysis is similar to \cite{Mairal_JMLR_2010,Bottou_CUP_1998}, but is different because of the update rule in (\ref{Eq:SVRMU_1_W}). More specifically, denoting the rightmost term in (\ref{Eq:problem}) as $l(\vec{h}_n, \mat{W}):=\frac{1}{2}\| \vec{v}_n - \mat{W} \vec{h}_n\|^2_2$, we define the {\it empirical cost} $f_N(\vec{h}_n,\mat{W})=\frac{1}{N}\sum_{n=1}^N l(\vec{h}_n, \mat{W})$. We also define $\hat{f}_N(\mat{W}):=\frac{1}{N}\sum_{n=1}^N l(\hat{\vec{h}}_n, \mat{W})$, where $\hat{\vec{h}}_n$ is already calculated during the previous steps. We now consider the {\it expected cost} $f(\vec{h}_n,\mat{W}):=\mathbb{E}_{\vec{v}}[l(\vec{h}_n, \mat{W})]= \lim_{N \rightarrow \infty} f_N(\vec{h}_n,\mat{W})$, where the expectation is taken with respect to the distribution $P(\vec{v})$ of the samples. Our interest is usually in the minimization of this expected cost $f(\vec{h}_n,\mat{W})$ almost surely (a.s.) instead of the empirical cost $f_N(\vec{h}_n,\mat{W})$. To this end, the convergence analysis first shows that $f_N(\vec{h}_n,\mat{W})-\hat{f}_N(\mat{W})$ converge a.s. to zero, where $\hat{f}_N(\mat{W})$ acts as a {\it surrogate function} for $f_N(\vec{h}_n,\mat{W})$. For this proof, we show that $\hat{f}_N(\mat{W})$ converges a.s. under the modified update rule in (\ref{Eq:SVRMU_1_W}). Here, the stepsize ratio $\alpha_t^s$ plays a crucial role in generating a diminishing sequence of $\mat{S}_t^s$ to guarantee its convergence. After showing the convergence of $f_N(\vec{h}_n,\mat{W})$, we finally obtain below;
\begin{Thm}
\label{Thm:GlobalConvAnalysis}
Assume that $\{\vec{v}\}_{n=1}^{\infty}$ are i.i.d. random processes, and bounded.
Iterates of $\mat{W}_{t}^{s}$ for $0\leq t\leq m_s-1$ and $0\leq s$ are compact. The initial $\tilde{\mat{W}}^0$ is nonnegative and has a full column rank. 
$\hat{f}_N(\mat{W})$ is positive definite and strictly convex.
$\alpha_t^s$ generates a diminishing stepsize of $\mat{S}_t^s$. 
Then, the iterates $\mat{W}_{t}^{s}$ produced by Algorithm \ref{Alg:SVRMU} asymptotically coincide with the stationary points of the minimization problem of $f(\vec{h}_n,\mat{W})$.
\end{Thm}

\section{Extensions of SVRMU}
This section proposes two variants of SVRMU.

\subsection{Accelerated SVRMU (SVRMU-ACC)}
\label{Sec:AccSVRMU}
Close examination of the update rule of $\vec{h}_{k}$ and $\mat{W}^s_{t}$ reveals that, whereas the latter requires $3FK+2FN$ because of the dominant calculation of the component-wise product of $\mat{W}^s_{t}$ at the last step, the former requires only $3FK+2K$, which is much lower than that of the latter because of $K \ll \{F, N\}$. Therefore, we can repeat the calculation of $\vec{h}_{k}$ several times, which corresponds to Step 7 in Algorithm \ref{Alg:SVRMU}, before the computation of $\mat{W}^s_{t}$. Although a similar strategy is also proposed for the batch-based MU \cite{Gillis_NN_2012}, the proposed one differs because of the different update rule. The noteworthy point is the stopping criteria, which are (i) the maximum iteration number $L$, and (ii) the dynamic stop criteria. The former specifically examines the ratio of the calculation complexity between $\mat{W}^s_{t}$ and $\vec{h}_{k}$. We calculate $L=\max\{\lfloor \beta\frac{3FK+2FN}{3FK+2K}\rfloor, 1\}$, where $0\leq \beta \leq 1$. Regarding the dynamic stop criteria, the process stops when the change between $l$-th $\vec{h}_{k}^{(l)}$ and $(l-1)$-th $\vec{h}_{k}^{(l-1)}$ falls below the predefined ratio $\epsilon$ of the difference from the initial value $\vec{h}_{k}^{(0)}$. The algorithm is presented as Algorithm \ref{Alg:SVRMU_ACC}.
\vspace*{-0.1cm}
\begin{algorithm}
\caption{Repetitive calculation algorithm of $\vec{h}_k$.}
\label{Alg:SVRMU_ACC}
\begin{algorithmic}[1]
\REQUIRE{$\vec{h}_k$, $(\mat{W}^s_{t})^T\vec{v}_k$, ${(\mat{W}^s_{t})}^T \mat{W}^s_{t}$ and the ratio $\epsilon$.}
\STATE{Set $\vec{h}_k^{(0)}=\vec{h}_k$.}
\FOR{$l=1,2,\ldots,L$} 
	\STATE{Calculate $\vec{h}_{k}=\vec{h}_{k} \odot {(\mat{W}^s_{t})}^T \vec{v}_k/({(\mat{W}^s_{t})}^T \mat{W}^s_{t} \vec{h}_{k})$. }
	\IF{$\| \vec{h}_k^{(l)} - \vec{h}_k^{(l-1)}\|_F < \epsilon \| \vec{h}_k^{(l)} - \vec{h}_k^{(0)} \|_F $}
		\STATE{break.}
	\ENDIF
\ENDFOR
\STATE{Return $\vec{h}_k=\vec{h}_k^{(l)}$.}
\end{algorithmic}
\end{algorithm}
\vspace*{-0.3cm}
\subsection{Robust SVRMU (R-SVRMU)}
\label{Sec:RobustVariants}
The outlier in \mat{V} causes remarkable degradation of the approximation of \mat{V}. To address this issue, the robust batch-NMF \cite{RobustManiNMF} and the robust online-NMF \cite{Zhao_ICASSP_2016} have been proposed. This extension also tackles the same problem within the SVRMU framework. 
Given the outlier matrix $\mat{R}=[\vec{r}_1, \ldots, \vec{r}_N] \in \mathbb{R}_+^{F \times N}$, the robust variant seeks $\mat{V}\approx \mat{WH} + \mat{R}$, of which minimization problem is formulated as 
\begin{eqnarray*}
	\label{Eq:problem_robust}
	\min_{\scriptsize \mat{W}, \mat{H}, \mat{R}} &&\frac{1}{N} \sum_{n=1}^N \frac{1}{2}\| \vec{v}_n - \mat{W} \vec{h}_n -\vec{r}_n\|^2_2 +  \lambda \| \vec{r}_n  \|_1,\nonumber \\
	 {\rm s.t.}&& [\mat{W}]_{f,k} \geq 0,\ \vec{h}_n \geq 0,\ \vec{r}_n \geq 0,\ \ \ \forall f,n,k,
\end{eqnarray*}
where $\lambda>0$ is the regularization parameter, and $\|\cdot \|_1$ is the $\ell{1}$-norm. For this problem, the update rule (\ref{Eq:SVRMU_1_W}) is redefined as
\begin{eqnarray*}
	\label{Eq:Robust_SVRMU}
	\mat{W}^s_{t+1}&=&\mat{W}^s_{k} - \mat{S}_t^s \odot \biggl[
	\left((\mat{W}^s_t\vec{h}_{k}+\vec{r}_{k})\vec{h}_{k}^{T} +  \vec{v}_{k} (\tilde{\vec{h}}_{k}^s)^{T} + \frac{(\tilde{\mat{W}}^s  \tilde{\mat{H}}^s + \tilde{\mat{R}}^s) (\tilde{\mat{H}}^s)^{T}}{N}\right) \nonumber \\
	&&- \left(\vec{v}_{k}\vec{h}_{k}^{T} + (\tilde{\mat{W}} \tilde{\vec{h}}^s_{k}+\tilde{\vec{r}}^s_k)(\tilde{\vec{h}}_{k}^s)^{T}+\frac{\mat{V} (\tilde{\mat{H}}^s)^{T}}{N} \right)
	\biggr].	\ \ \ \ \ \ \ \ 
\end{eqnarray*}
Accordingly, we respectively calculate as
\begin{eqnarray*}
 \vec{h}_k &\leftarrow&\vec{h}_{k}\odot\frac{(\mat{W}_{t}^s)^T \vec{v}_{k}}{(\mat{W}_{t}^s)^T \mat{W}^s_{t} \vec{h}_{k} + (\mat{W}^s_{t})^T \vec{r}_t}\\
 \vec{r}_{k}  &\leftarrow&  \mat{r}_{k} \odot \frac{\vec{v}_{k}}{\mat{W}_{t}^s \vec{h}_{k} + \vec{r}_k + {\bf \Lambda}_{F\times 1}}. 
\end{eqnarray*}
\section{Numerical experiments}
\label{Sec:Numerical_experiments}

This section demonstrates the effectiveness of SVRMU by comparing with the state-of-the-art online algorithms for NMF. 
We implemented all of the algorithms in Matlab\footnote{\url{https://github.com/hiroyuki-kasai}}. 

\subsection{Convergence behavior under clear synthetic data}
\begin{figure}[htbp]
\vspace*{-0.3cm}
\begin{center}
\begin{tabular}{c}
\begin{minipage}{0.5\hsize}
	\begin{center}
	\hspace*{-0.5cm}\includegraphics[width=\hsize]{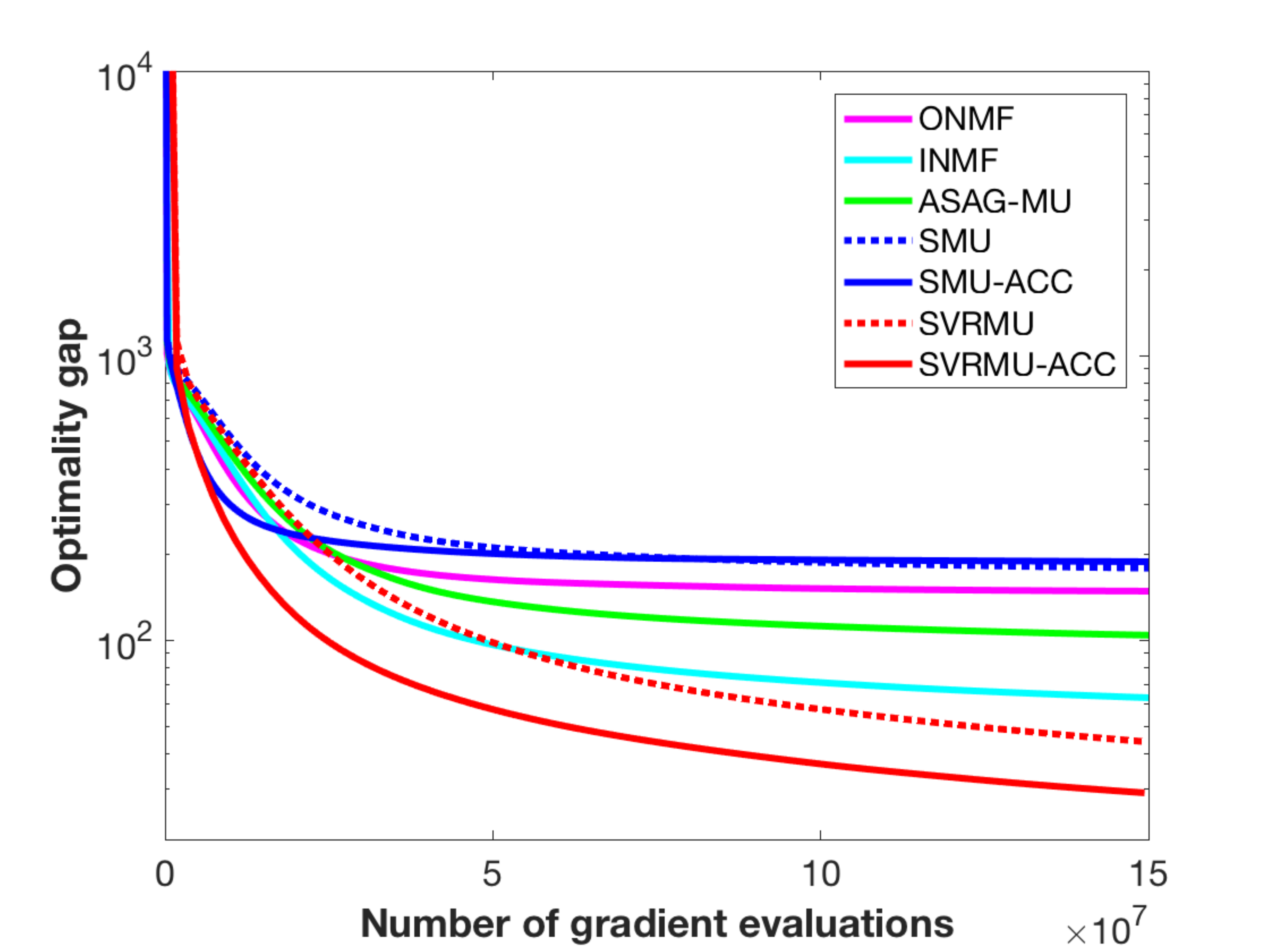}
	
	{\scriptsize (a) \# of gradient counts v.s. optimality gap}
	\end{center}
\end{minipage}
\begin{minipage}{0.5\hsize}
	\begin{center}
	\hspace*{-0.3cm}\includegraphics[width=\hsize]{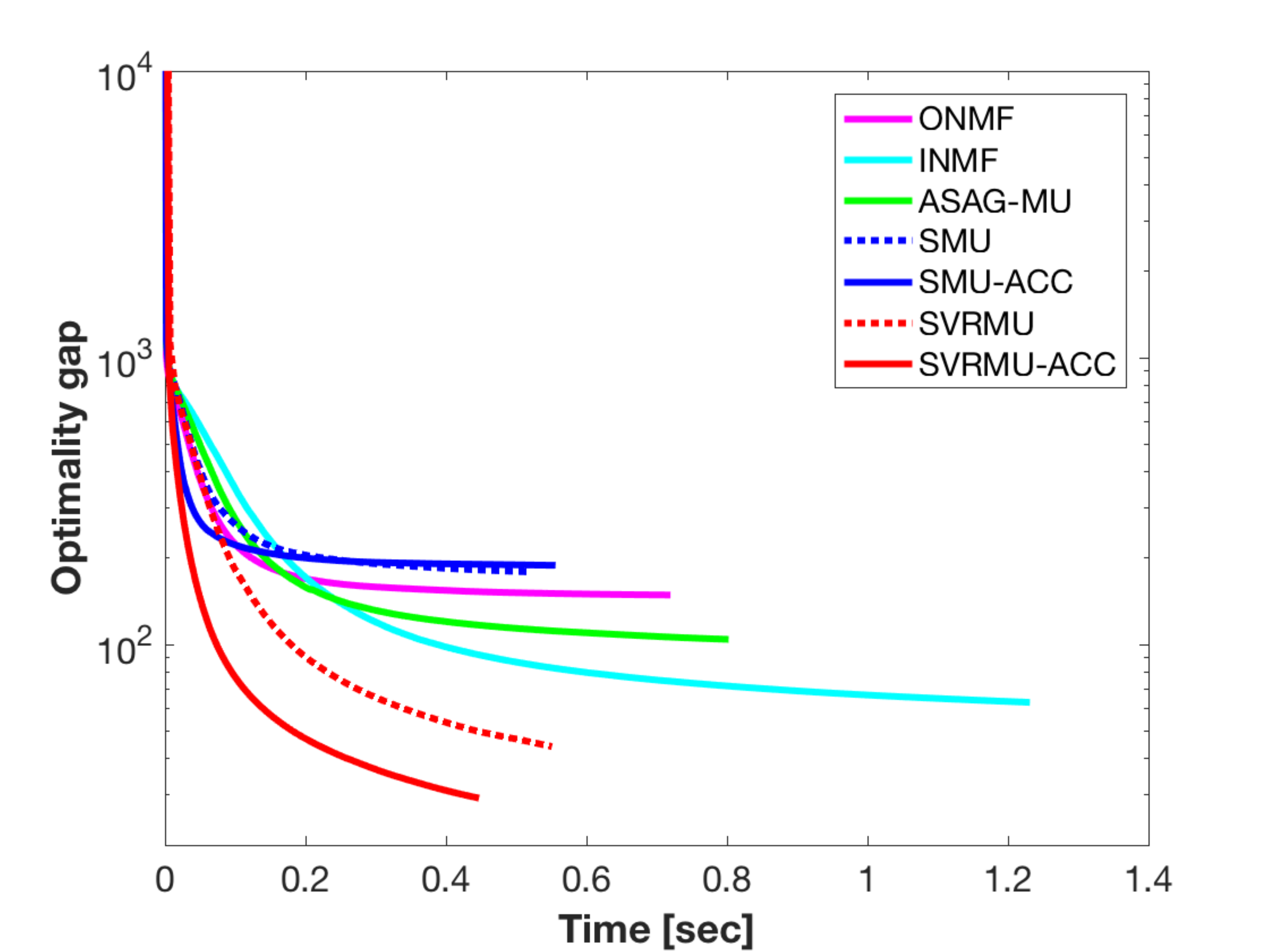}
	
	{\scriptsize (b) Time v.s. optimality gap (enlarged)}	
	\end{center}	
\end{minipage}
\end{tabular}
\end{center}
\caption{Convergence behavior on synthetic dataset.}
\label{fig:synthetic}
\end{figure}

The element $[\mat{W}_o]_{f,n}$ of the ground-truth $\mat{W}_o \in \mathbb{R}^{F \times K_o}_+$  is generated from a Gaussian distribution with a mean of zero and  variance $1/\sqrt{K_o}$ for any $(f,n)$, where $K_o$ is the ground-truth rank dimension. Similarly, we generate $\mat{H}_o \in \mathbb{R}^{K_o \times N}_+$. Then, the clean data $\mat{V}_o$ are created as $\mat{V}_o=\mathcal{P}_{\tilde{\mathcal{V}}}(\mat{W}_o \mat{V}_o)$, where $\mathcal{V}=[0,1]^{F \times N}$, and where $\mathcal{P}_{\tilde{\mathcal{V}}}$ is the normalization projector \cite{Zhao_ICASSP_2016}. We set $(F,N,K_o,b)=(300,1000,10,100)$. The maximum epoch is $500$.
The following methods are used for comparison: incremental MU (INMF) \cite{Bucak_PR_2009}, online MU (ONMF) \cite{Zhao_ICASSP_2016}, and ASAG-MU \cite{Serizel_MLSP_2016}. 
Our proposed algorithms include 
SMU and SVRMU in Section \ref{Sec:SVRMU}, and those accelerated variants, i.e., SMU-ACC and SVRMU-ACC, in Section \ref{Sec:AccSVRMU}.
Figure \ref{fig:synthetic} presents results of the convergence behavior in terms of the {\it optimality gap}, which is calculated using HALS \cite{Cichocki_IEICETranA_2009} in advance. The figure shows the superior performance of SVRMU in terms of the number of gradients and the time. 

\subsection{Base representation of face image with outlier}
\begin{figure}[t]
\begin{center}
\begin{tabular}{c}
\hspace*{-0.5cm}\begin{minipage}{0.4\hsize}
	\begin{center}
	\includegraphics[width=\hsize]{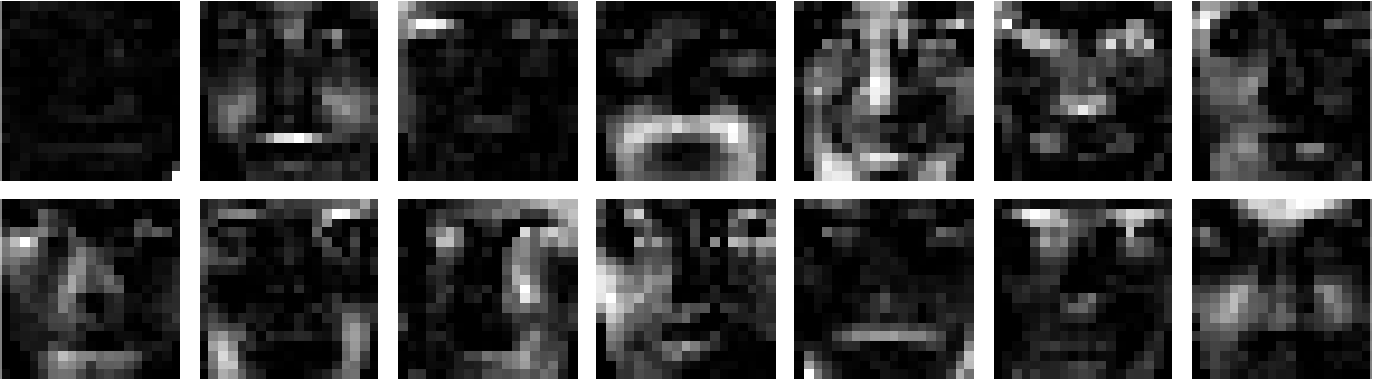}	
	\vspace*{-0.6cm}
	
	{\small \hspace*{0.6cm} (a) R-NMF (batch)}
	\end{center}
\end{minipage}
\hspace*{0.9cm}
\begin{minipage}{0.4\hsize}
	\begin{center}
	\includegraphics[width=\hsize]{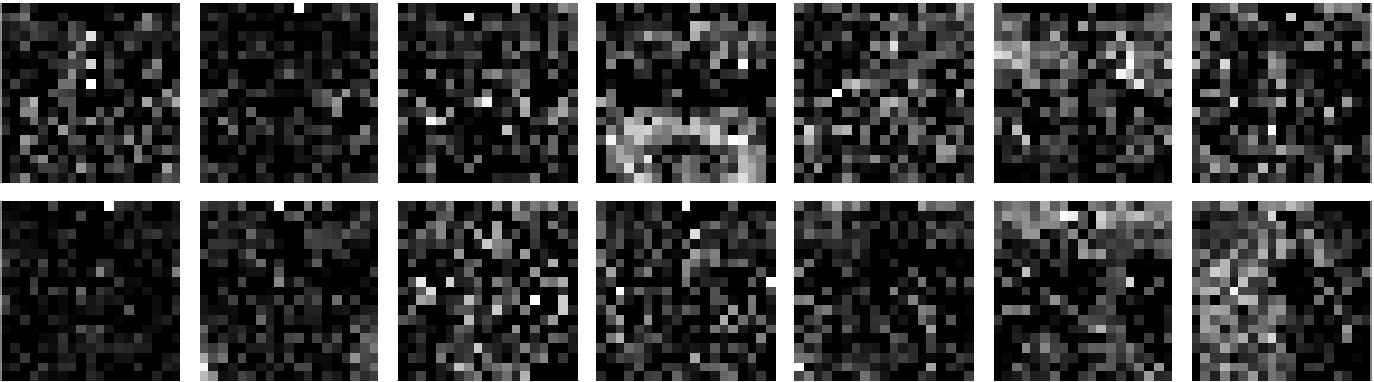}		
	\vspace*{-0.6cm}
	
	{\small \ \ \ (b) ONMF}	
	\end{center}	
\end{minipage}\\
\vspace*{-0.3cm}\\

\hspace*{-0.5cm}\begin{minipage}{0.4\hsize}
	\begin{center}
	\includegraphics[width=\hsize]{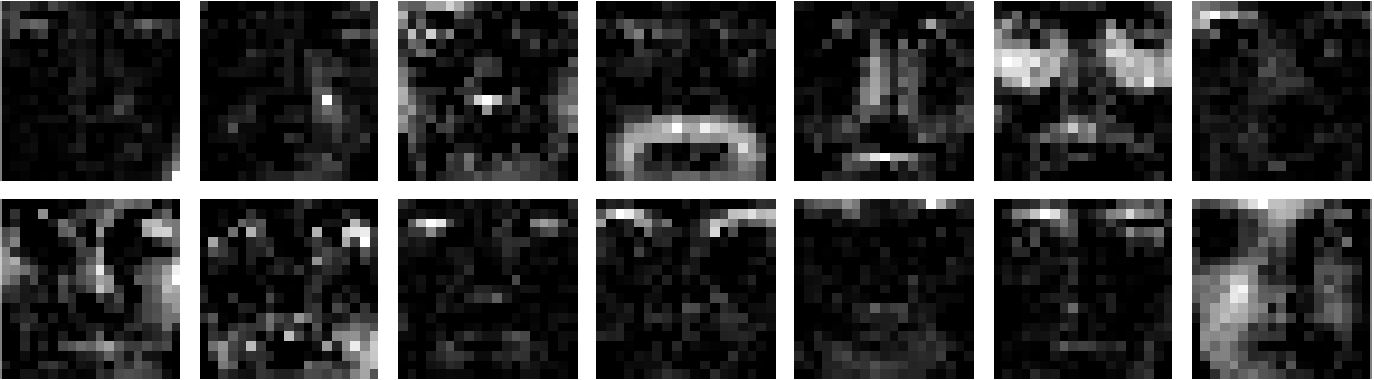}		
	\vspace*{-0.6cm}	
	
	{\small \hspace*{0.6cm} (c) R-ONMF}		
	\end{center}
\end{minipage}
\hspace*{0.9cm}
\begin{minipage}{0.4\hsize}
	\begin{center}
	\includegraphics[width=\hsize]{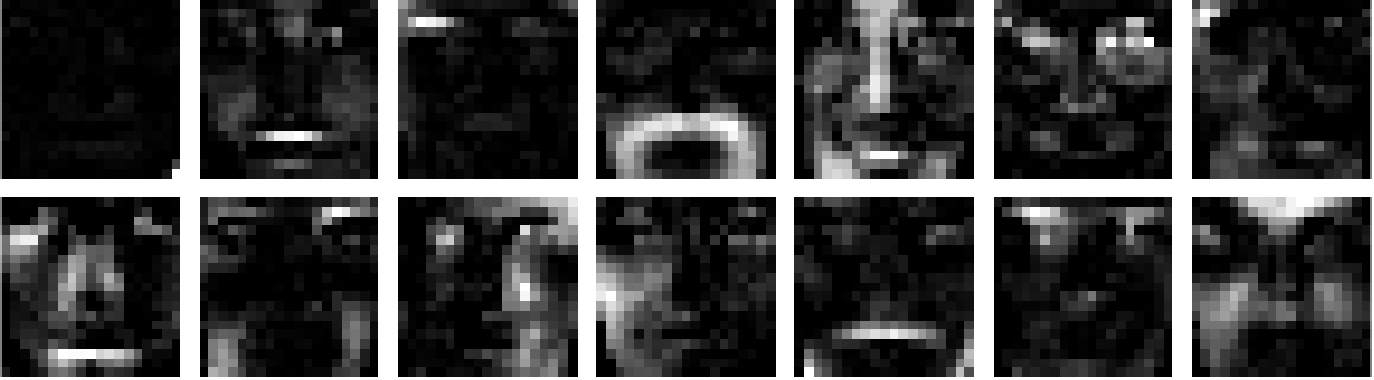}			
	\vspace*{-0.6cm}	
	
	{\small (d) R-SVRMU (proposed)}		
	\end{center}	
\end{minipage}
\end{tabular}
\caption{Basis representations on the CBCL dataset.}
\label{fig:CBCL_BaseRep}
\end{center}

\end{figure}

We use the CBCL face dataset\footnote{\url{http://cbcl.mit.edu/cbcl/software-datasets/FaceData2.html}}, which has $2429$ gray-scale images of size $19\times 19$. The maximum level of the pixel values is set to 50. All pixel values are normalized. We also randomly add entry-wise nonnegative outliers with density $\rho=0.9$. All outliers are drawn from the i.i.d. from a uniform distribution $\mathcal{U}[30,50]$. $K_o$ is fixed to $49$.
The methods of comparison include ONMF and its robust variant: R-ONMF \cite{Zhao_ICASSP_2016}, and the accelerated variant of R-SVRMU. The batch-based variant of R-ONMF (R-NMF) is also evaluated. Figure \ref{fig:CBCL_BaseRep} presents an illustration of the generated 14 basis representations, where it is apparent that R-SVRMU produces better bases than R-ONMF, and gives similar bases as the batch-based R-NMF.

\vspace*{1.0cm}

\section{Conclusions}
\label{Sec:Conclusions}
This present paper has proposed a novel stochastic multiplicative update with variance reduction technique: SVRMU. Numerical comparisons suggest that SVRMU robustly outperforms state-of-the-art algorithms across different synthetic and real-world datasets.

\clearpage
\bibliographystyle{unsrt}
\bibliography{/Users/kasai/Dropbox/DOC/Research/bibtex/nmf,/Users/kasai/Dropbox/DOC/Research/bibtex/stochastic_online_learning}

\begin{thebibliography}{10}

\bibitem{Lee2001}
Daniel~D Lee and H.~Sebastian Seung.
\newblock Algorithms for non-negative matrix factorization.
\newblock {\em Advances in neural information processing systems (NIPS)}, 2001.

\bibitem{Lin_IEEENN_2007}
C.-J. Lin.
\newblock On the convergence of multiplicative update algorithms for
  nonnegative matrix factorization.
\newblock {\em IEEE Transactions on Neural Networks}, 18(6):1589--1596, 2007.

\bibitem{Hibi_ICONIP_2011}
R.~Hibi and N.~Takahashi.
\newblock A modified multiplicative update algorithm for euclidean
  distance-based nonnegative matrix factorization and its global convergence.
\newblock In {\em International Conference on Neural Information Processing
  (ICONIP)}, pages 655--662, 2011.

\bibitem{Lin_NC_2007}
C.-J. Lin.
\newblock Projected gradient methods for non-negative matrix factorization.
\newblock {\em Neural Comput.}, 19(10):2756--2779, 2007.

\bibitem{Cichocki_IEICETranA_2009}
A.~Cichocki and P.~Anh-Huy.
\newblock Fast local algorithms for large scale nonnegative matrix and tensor
  factorizations.
\newblock {\em IEICE Transactions on Fundamentals of Electronics,
  Communications and Computer Sciences}, 92(3):708--721, 2009.

\bibitem{Gillis_NN_2012}
N.~Gillis and F.~Glineur.
\newblock Accelerated multiplicative updates and hierarchical als algorithms
  for nonnegative matrix factorization.
\newblock {\em Neural Comput.}, 24(4):1085--1105, 2012.

\bibitem{Kim_JGO_2014}
J.~Kim, Y.~He, and H.~Park.
\newblock Algorithms for nonnegative matrix and tensor factorizations: A
  unified view based on block coordinate descent framework.
\newblock {\em Journal of Global Optimization}, 58(2):285--319, 2014.

\bibitem{Bucak_PR_2009}
S.~S. Bucak and B.~Gunsel.
\newblock Incremental subspace learning via non-negative matrix factorization.
\newblock {\em Pattern Recognition}, 42(5):788--797, 2009.

\bibitem{Fvotte_NC_2009}
C.~F{\'e}votte, N.~Bertin, and JL~Durrieu.
\newblock Nonnegative matrix factorization with the itakura-saito divergence:
  with application to music analysis.
\newblock {\em Neural Comput.}, 21(3):793--830, 2009.

\bibitem{Guan_IEEENNLS_2012}
N.~Guan, D.~Tao, Z.~Luo, and B.~Yuan.
\newblock Online nonnegative matrix factorization with robust stochastic
  approximation.
\newblock {\em IEEE Trans. Neural Netw. Learn. Syst.}, 23(7):1087, 1099 2012.

\bibitem{Zhao_AISTATS_2017}
R.~Zhao, V.~Y.~F. Tan, and H.~Xu.
\newblock Online nonnegative matrix factorization with general divergences.
\newblock In {\em International Conference on Artificial Intelligence and
  Statistics (AISTATS)}, 2017.

\bibitem{Zhao_ICASSP_2016}
R.~Zhao and V.~Y.~F. Tan.
\newblock Online nonnegative matrix factorization with outliers.
\newblock In {\em IEEE International Conference on Acoustics, Speech, and
  Signal Processing (ICASSP)}, 2016.

\bibitem{Robbins_MathStat_1951}
H.~Robbins and S.~Monro.
\newblock A stochastic approximation method.
\newblock {\em Ann. Math. Statistics}, pages 400--407, 1951.

\bibitem{Johnson_NIPS_2013_s}
R.~Johnson and T.~Zhang.
\newblock Accelerating stochastic gradient descent using predictive variance
  reduction.
\newblock In {\em NIPS}, pages 315--323, 2013.

\bibitem{Roux_NIPS_2012_s}
N.~L. Roux, M.~Schmidt, and F.~R. Bach.
\newblock A stochastic gradient method with an exponential convergence rate for
  finite training sets.
\newblock In {\em NIPS}, pages 2663--2671, 2012.

\bibitem{Shalev_JMLR_2013_s}
S.~Shalev-Shwartz and T.~Zhang.
\newblock Stochastic dual coordinate ascent methods for regularized loss
  minimization.
\newblock {\em JMLR}, 14:567--599, 2013.

\bibitem{Defazio_NIPS_2014_s}
A.~Defazio, F.~Bach, and S.~Lacoste-Julien.
\newblock {SAGA}: A fast incremental gradient method with support for
  non-strongly convex composite objectives.
\newblock In {\em NIPS}, 2014.

\bibitem{Zhang_SIAMJO_2014_s}
Y.~Zhang and L~Xiao.
\newblock Stochastic primal-dual coordinate method for regularized empirical
  risk minimization.
\newblock {\em SIAM J. Optim.}, 24(4):2057--2075, 2014.

\bibitem{Nguyen_ICML_2017}
L.~M. Nguyen, J.~Liu, K.~Scheinberg, and M.~Takac.
\newblock {SARAH}: A novel method for machine learning problems using
  stochastic recursive gradient.
\newblock In {\em ICML}, 2017.

\bibitem{Johnson_NIPS_2013}
R.~Johnson and T.~Zhang.
\newblock Accelerating stochastic gradient descent using predictive variance
  reduction.
\newblock In {\em Advances in Neural Information Processing Systems (NIPS)},
  pages 315--323, 2013.

\bibitem{Mairal_JMLR_2010}
J.~Mairal, F.~Bach, J.~Ponce, and G.~Sapiro.
\newblock Online learning for matrix factorization and sparse coding.
\newblock {\em Journal of Machine Learning Research (JMLR)}, 11:19--60, 2010.

\bibitem{Bottou_CUP_1998}
L.~Bottou.
\newblock Online algorithm and stochastic approximations.
\newblock In David Saad, editor, {\em On-Line Learning in Neural Networks}.
  Cambridge University Press, 1998.

\bibitem{RobustManiNMF}
Huang Jin, NIE Feiping, and Ding Chris.
\newblock Robust manifold nonnegative matrix factrization.
\newblock {\em ACM Transations on Knowledge Discovery from Data}, 8(3):1--21,
  2014.

\bibitem{Serizel_MLSP_2016}
Romain Serizel, Slim Essid, and Ga{\"e}l Richard.
\newblock Mini-batch stochastic approaches for accelerated multiplicative
  updates in nonnegative matrix factorisation with beta-divergence.
\newblock In {\em IEEE International Workshop on Machine Learning for Signal
  Processing (MLSP)}, pages 5470--5474, 2016.

\end{thebibliography}
\end{document}